\documentclass{svproc}

\usepackage{soul}
\usepackage{amsmath}
\usepackage{amssymb}
\usepackage{bm}
\usepackage{tabu}
\usepackage{xcolor}
\usepackage{hyperref}
\usepackage{tikz}
\usepackage{newtxtext,newtxmath}
\usepackage{caption}
\usepackage{subcaption}
\usepackage{soul}
\usepackage[colorinlistoftodos]{todonotes}

\usepackage{url}

\begin{document}
\mainmatter              % start of a contribution
\newcommand{\bbmu}{\bm{\mu}}
\newcommand{\bbz}{\bm{U}}
\newcommand{\bbu}{\bm{u}}
\newcommand{\bbv}{\bm{v}}
\newcommand{\bbd}{\bm{d}}
\newcommand{\monica}[1]{\textcolor{magenta}{#1}}
\newcommand{\ivan}[1]{\textcolor{orange}{#1}}

\title{A time-adaptive algorithm for pressure dominated flows: a heuristic estimator}
\titlerunning{A time-adaptive algorithm for pressure dominated flows: a heuristic estimator}  % abbreviated title (for running head)
\author{Ivan Prusak\inst{1} \and Davide Torlo\inst{2} \and
Monica Nonino\inst{3} \and Gianluigi Rozza \inst{4} }
\authorrunning{Ivan Prusak et al.} % abbreviated author list (for running head)
\tocauthor{Ivan Prusak, Davide Torlo, Monica Nonino, Gianluigi Rozza}
\institute{Sissa, Mathematics Area, mathLab, International School for Advanced Studies, \\ via Bonomea 265, 34136 Trieste, Italy, \email{iprusak@sissa.it}
\and
Dipartimento di Matematica ``Guido Castelnuovo'', Universit\'a di Roma La Sapienza, piazzale Aldo Moro 5, 00185 Roma, Italy, \email{davide.torlo@uniroma1.it}
\and
Faculty of Mathematics, University of Vienna, Oskar-Morgenstern-Platz 1, \\1090, Vienna, Austria, \email{monica.nonino@univie.ac.at}
\and Sissa, Mathematics Area, mathLab, International School for Advanced Studies, \\ via Bonomea 265, 34136 Trieste, Italy, \email{grozza@sissa.it}}

\maketitle              % typeset the title of the contribution
\begin{abstract}
% The abstract should summarize the contents of the paper
% using at least 70 and at most 150 words. It will be set in 9-point
% font size and be inset 1.0 cm from the right and left margins.
% There will be two blank lines before and after the Abstract. \dots
% We would like to encourage you to list your keywords within
% the abstract section using the \keywords{...} command.
% \st{Time-adaptive stepping schemes are indispensable in applications that require accurate solutions of complex physical phenomena over extended periods.} 
This work aims to introduce a heuristic timestep-adaptive algorithm for Computational Fluid Dynamics (CFD) and Fluid-Structure Interaction (FSI) problems where the flow is dominated by the pressure. In such scenarios, many time-adaptive algorithms based on the interplay of implicit and explicit time schemes fail to capture the fast transient dynamics of pressure fields. 
We present an algorithm that relies on a temporal error estimator using Backward Differentiation Formulae (BDF$k$) of order $k=2,3$. Specifically, we demonstrate that the implicit BDF$3$ solution can be well approximated by applying a single Newton-type nonlinear solver correction to the implicit BDF$2$ solution. The difference between these solutions determines our adaptive temporal error estimator. The effectiveness of our approach is confirmed by numerical experiments conducted on a backward-facing step flow CFD test case with Reynolds number $300$ and on a two-dimensional haemodynamics FSI benchmark.  

\end{abstract}

\section{Introduction}
% Descrizione problema (incluso FEM NS + FEM FSI), necessità di time adaptivity and literature for time adapt (in general) and for FSI-NS
% \begin{itemize}
%     \item \st{CFD, FSI}
%     \item \st{Why time adaptive}
%     \item Literature: fractional timesteps, explicit-implicit (Adaptive time-stepping for FSI, M. Mayra, W. A. Wallc, M. W. Geeb)
%     \item \st{Problems: more solutions, not get pressure}
%     \item \st{Outline}
% \end{itemize}
In this manuscript, we are interested in solving nonlinear problems that arise from computational fluid dynamics (CFD) and fluid-structure interaction (FSI) problems. We denote by $\Omega\subset\mathbb{R}^2$ the physical domain of interest; we further assume that $\partial\Omega=\Gamma_D\cup\Gamma_N$, where $\Gamma_D$ is the Dirichlet boundary and $\Gamma_N$ is the Neumann one. Let $[0,T]$, $T>0$ be the time interval. For every $t\in(0,T]$, we seek the solution $\bbz(t)\in W$ to the following PDE
\begin{equation}\label{eq:general_nonlinear_problem}
    R(\partial_t\bbz(t), \bbz(t))=0 \quad \text{in} \ W_0',
\end{equation}
where $R$ is the residual operator of our nonlinear PDE and $W$ and $W_0$ are appropriate trial and test function spaces, respectively. We can now introduce the weak formulation of problem~\eqref{eq:general_nonlinear_problem} corresponding to a CFD and a FSI problem, respectively.\\ %Since we are working in a nonlinear setting, the weak formulation will be presented in terms of the residual function.\\
\textbf{Navier-Stokes problem.}
Let $\bbu$ and $p$ be the two components of the solution $\bbz$ introduced in Eq.~\eqref{eq:general_nonlinear_problem}, namely the fluid velocity and the fluid pressure. We introduce the following function spaces $W = V\times Q$ and $W_0 = V_0\times Q$ with
\begin{equation*}
    \begin{split}
        V &:= \{\bbv\in [H^1(\Omega)]^2 \text{ s.t. }\bbv=\bbu_D\text{ on }\Gamma_D\},\\
        V_0 &:= \{\bbv\in [H^1(\Omega)]^2 \text{ s.t. }\bbv=\bm{0}\text{ on }\Gamma_D\},\\
        Q &:= L^2(\Omega),
    \end{split}
\end{equation*}
where $\bbu_D$ is a given Dirichlet datum on $\Gamma_D \times (0,T].$ %\subseteq \partial \Omega$.
The weak formulation of the Navier-Stokes problem reads as follows: for every $t\in (0, T]$, find $\bbu(t)\in V$ and $p(t)\in Q$ such that, for every $\bbv\in V_0$ and every $q\in Q$ the following holds
\begin{equation}\label{eq:weak_time_continuous_ns}
    \begin{cases}
        m(\partial_t\bbu, \bbv)+a(\bbu, \bbv)+c(\bbu, \bbu, \bbv) + b(p,\bbv)=\left( \bbu_N, \bbv\right)_{\Gamma_N} \quad\forall t\in (0, T],\\
        b(\bbu, q)=0 \quad\forall t\in (0, T],\\
        \bbu(t=0)=\bbu_0\quad\text{in }\Omega,
    \end{cases}
\end{equation}
where $\bbu_0$ is a given initial condition (IC) and $\bbu_N$ is a Neumann datum defined on $\Gamma_N\times(0,T]$.
In system~\eqref{eq:weak_time_continuous_ns}, $(\cdot, \cdot)_w$ indicates $L^2(\omega)$ inner product and the bilinear and trilinear forms are introduced in Table~\ref{tab:ns_forms}. 
\begin{table}
\caption{Navier-Stokes forms definitions}   \label{tab:ns_forms}
\begin{center}
\bgroup
\def\arraystretch{1.5}
\begin{tabular}{|c|c|}\hline
% \centering
%    \tabulinesep=0.5mm
%     \begin{tabu}{|c|l|}\hline
        Form & Definition \\    \hline 
        fluid velocity mass &  $m(\bm{w},\bbv) = \int_{\Omega} \bm{w} \cdot \bbv d\Omega$\\ \hline 
        fluid stiffness & $a(\bbu,\bbv) = \int_{\Omega} \nu\nabla\bbu:\nabla\bbv d\Omega$ \\ \hline 
        fluid incompressibility & $b(p, \bbv) = -\int_{\Omega}p\text{ div}\bbv d\Omega$ \\ \hline   
        fluid advection &$ c(\bbu, \bm{w},\bbv) = \int_{\Omega} \left( \bbu \cdot\nabla \right)  \bm{w} \cdot \bbv d\Omega$ \\ \hline
        %fluid forcing & $f_f( v_f; d_f) = \int_{\Omega^f} J b_f \cdot v_f d\Omega$\\ \hline
    % \end{tabu}
\end{tabular}
\egroup
\end{center}
\end{table}

\textbf{Fluid-Structure Interaction problem.}
Due to the different nature of the formalism used to describe the motion of a fluid (Eulerian formalism) and the motion of a solid (Lagrangian formalism), in this manuscript we make use of the so-called ALE formalism; for more details on this, we refer to~\cite{richter2017,prusakPhDthesis}. Hereafter, we will assume that all the variables are considered on a reference configuration $\Omega$, which is time-independent. We will resort to the monolithic Lagrangian formulation as described in~\cite{ballarin2016pod}. 
For an FSI problem, we assume a partition of $\Omega$ into two disjoint open domains $\Omega^f$ (the fluid subdomain) and $\Omega^s$ (the solid subdomain) and we denote by $\Gamma_I:=\Bar{\Omega}^f\cap\Bar{\Omega}^s$ the fluid-structure interface. 
We further denote with $\Gamma_{f,D}$ and $\Gamma_{s,D}$ the Dirichlet boundaries, and with $\Gamma_{f,N}$ and $\Gamma_{s,N}$ the Neumann boundaries, for the fluid and the solid, respectively. 
The components of $\bbz$, for an FSI problem, are the following: the common velocity field for the fluid and the structure subproblems $\bbu = (\bbu_f, \bbu_s)\colon\Omega\rightarrow \mathbb{R}^2$, the fluid pressure $p_f\colon\Omega^f\rightarrow\mathbb{R}$ and the common displacement field $\bbd = (\bbd_f, \bbd_s)\colon\Omega\rightarrow\mathbb{R}^2$. We introduce the following function spaces $W = W_0 = V\times Q \times E$:
\begin{eqnarray*}
    V & := & \left\{\bbv \in \left[ H^1(\Omega)\right]^2: \bbv = 0 \text{ on } \Gamma_{f,D}  \ \text{and $\bbv$ is continuous across } \Gamma_I   \right\}, \\
    Q & := & L^2(\Omega^f), \\
     E& := &\left\{\bm{e}  \in \left[ H^1(\Omega)\right]^2: \bm{e} = 0 \text{ on } \Gamma_{D} \cup \Gamma_{f,N}   \ \text{and $\bm{e}$ is continuous across } \Gamma_I   \right\}.
\end{eqnarray*}

The weak formulation of the coupled FSI problem reads as follows: for every $t\in(0, T]$, find $(\bbu(t), p_f(t), \bm{d}(t)) = (\bbu_f(t),\bbu_s(t), p_f(t), \bm{d}_f(t), \bm{d}_s(t))\in V \times Q \times E$, such that the following holds: for every $\bbv = (\bbv_f, \bbv_s) \in V$, $q_f\in Q$ and $\bm{e} = (\bm{e}_f, \bm{e}_s)\in E$ 
\begin{equation}\label{eq:weak_time_continuous_fsi}
    \begin{cases}
        m_f(\partial_t\bbu, \bbv; \bbd) + m_s(\partial_t\bbu, \bbv) + a_f(\bbu, \bbv; \bbd)  +   a_s(\bbd, \bbv) + b_f^A(p_f, \bbv; \bbd) \\ \quad  + c_f^{ALE}(\partial_t\bbd, \bbv, \bbu; \bbd)          +  c_f(\bbu, \bbu, \bbv; \bbd)  = (\bbu_N^f, \bbv_f)_{\Gamma_{f, N}} +  \left(\bbd_N^s, \bbv_s\right)_{\Gamma_{s, N}},\\
        b_f^B(\bbu, q_f; \bbd) = 0,\\
        a_f^e(\bbd, \bm{e})  + a_f^{e, \partial f}(\bbd, \bm{e}) + \left(\partial_t \bbd_s - \bbu_s, \bm{e}_s\right)_{\Omega^s}   = 0,
    \end{cases}
\end{equation}
where $(\cdot, \cdot)_\omega$ indicates $L^2(\omega)$ inner product and all forms are defined in Table~\ref{tab:fsi_forms}.
\begin{table}
\caption{FSI forms definitions}\label{tab:fsi_forms} 
\begin{center}
\bgroup
\def\arraystretch{1.5}
\begin{tabular}{|c|c|}\hline
    % \begin{tabu}{|c|l|}\hline
        Form & Definition \\    \hline 
        fluid velocity mass &  $  m_f(\bm{w}, \bbv; \bbd) = \int_{\Omega^f} J\rho_f \bm{w} \cdot \bbv d\Omega$\\ \hline 
        fluid stiffness & $  a_f(\bbu, \bbv; \bbd) = \int_{\Omega^f} J\sigma_f^{du}(\bbu_f; \bbd_f) F^{-T} : \nabla \bbv_f d\Omega$ \\ \hline 
        fluid ALE & $  c_f^{ALE}(\bm{w}, \bbv, \bbu_f;\bbd) = -\int_{\Omega^f} J\rho_f \left[ \nabla \bbu_f F^{-1}\right]\bm{w} \cdot \bbv_f d\Omega$ \\ \hline
        fluid incompressibility A & $  b_f^A(p_f, \bbv; \bbd) = \int_{\Omega^f} J\sigma_f^{p}(p_f) F^{-T} : \nabla \bbv_f  d\Omega$ \\ \hline   
        fluid incompressibility B & $ b_f^B(\bbu, q_f; \bbd) = -\int_{\Omega^f} \text{div}\left( JF^{-1}\bbu_f\right) q_f d\Omega$ \\ \hline
        fluid advection &$ c_f(\bbu_f, \bm{w}_f, \bbv_f; \bbd_f) = \int_{\Omega^f} J\rho_f \left[ \nabla \bm{w}_f F^{-1}\right]\bbu_f \cdot \bbv_f d\Omega$ \\ \hline
        extension fluid stiffness& $  a_f^e(\bbd, \bm{e}) = \int_{\Omega^f} \nabla \bbd_f : \nabla \bm{e}_f d\Omega$\\ \hline
        extension interface stiffness& $  a_f^{e, \partial f}(\bbd, \bm{e}) = - \int_{\Gamma_I} \left[ \nabla \bbd_f \right] \bm{n}_f:  \bm{e}_f dS$\\ \hline
        structure displacement mass& $ m_s(\bbd, \bm{e}) = \int_{\Omega^s} \rho_s \bbd_s\cdot \bm{e}_s d\Omega$\\ \hline
        structure stiffness& $ a_s(\bbd, \bbv) = \int_{\Omega^s} P(\bbd_s) : \nabla \bbv_s d\Omega$\\ \hline
    % \end{tabu}
\end{tabular}
\egroup
\end{center}
\end{table}

There, $\rho_f$ and $\rho_s$ represent the fluid and the solid density, respectively.
System~\eqref{eq:weak_time_continuous_fsi} is then completed by some suitable ICs on $\bbu_f$, $\bbd_s$ and $\partial_t\bbd_s$.
In system~\eqref{eq:weak_time_continuous_fsi}, $\bbu_N^f$ and $\bbd_N^s$ are given Neumann data for the fluid and the solid; $\sigma_f(\bbu_f, p_f;\bbd_f)=\sigma_f^{du}(\bbu_f;\bbd_f)+\sigma_f^p(p_f)$ is the fluid Cauchy stress tensor, $P(\bbd_s)$ is the Piola-Kirchhoff tensor for the solid, $F$ is the Jacobian of the ALE map and $J$ is the determinant of $F$ (see \cite{richter2017} for more details); finally, $\bm{n}_f$ and $\bm{n}_s$ are the normals to $\Gamma_I$, outgoing the fluid and the solid subdomain, respectively.  

At this stage, we presented the problems of interest in the time-continuous setting: no time-stepping scheme has been introduced so far. 

\textbf{Overview of time--adaptive techniques.}
The choice of the time-step in numerical solvers for ODEs and PDEs has always been a crucial topic. Already from 1911, Richardson \cite{richardson1911ix} proposed an extrapolation procedure to estimate the error during a single time-step and to use it to correct its length.
This procedure was based on a second simulation run within the same time-step but using two time-steps of half the size of the original one. 
The resulting difference between the two simulations is of the order of the error of the simulation with the large time-step and, hence, gives an estimation of the local error.
Many estimators are still based on this idea.
A generalization of this concept involves to two different simulations in one time-step with different accuracies \cite{merson1957operational}. From there on, various \textit{embedded methods} were developed. 
They exploit the stages/steps structure to create two different approximations, typically of different orders.  Through their difference it is possible to obtain an error estimate for the worst approximated solution \cite{hairer1987solving}.
Commonly, a tolerance is set for the local error and the time-step is adjusted to fulfil this error bound in the current or following time-step. Various strategies have been used to set the time-step and other regularity constraints may be enforced not to vary the time-step too much \cite{hairer1987solving}.
When dealing with saddle-point simulations, most of these strategies fail at bounding errors for all the variables, such as the pressure in FSI and NS.
Hence, different options are available and the choice of the error estimator leads to great changes in the simulations.

In Section~\ref{sec:time_estimator}, we will describe a couple of error estimators based on the difference between BDF2 and BDF3 formulations which is able to address the issues stated above and provide an extensive numerical analysis of the proposed techniques on various CFD and FSI benchmark problems in Section~\ref{sec:simulations}.

\section{Time discretization and error estimator}\label{sec:time_estimator}
At the beginning of this section, we will briefly describe Backward Differentiation Formulae (BDF$k$) of order $k=2,3$ for non--constant time discretisation. Then, we will present a temporal error estimator based on the interplay of BDF$2$ and BDF$3$ implicit schemes and discuss the efficient approximation of the implicit error estimator.
% \todo{iv: da riscrivere nel modo piu' bello } 
% First of all, let us introduce the BDF formulation for the FSI and NS problems.

% \begin{itemize}
    % \item \st{BDF2 BDF3 , variable timesteps}
    % \item con focus su pressione che manca in altri stimatori
    % \item Algorithm (Adaptive time-stepping for FSI, M. Mayra, W. A. Wallc, M. W. Geeb) Formula di $\Delta t$, parametri.
    % \item Domanda: quale error estimator?
    % \item stimatore esplicito ha problemi con grandi differenze in pressione
%     \item BDF3-BDF2 implicit buono ma costoso, dove possiamo risparmiare? Matrix assembly può essere risparmiato? 
%     \item Semplifichiamo Newton di BDF3 con la matrice già calcolata per fare 1 iterazione di Newton per avere stesso errore
%     \item proof with Taylor expansion (per NS)
% \end{itemize}

\subsection{Backward Differentiation Formulae (BDF) with varying time-steps} \label{sec:bdf_k}
We assume the following partition of the time interval $[0,T]$: $0=:t_0 < t_1 < \dots < t_M :=T, M>0$ and denote by $\Delta_t^n:= t_n - t_{n-1}$ for $n=1,\dots,M$. Furthermore, we introduce the following notation: $\bbz^n := \bbz(t_n), n=0,\dots, M$. In the BDF, the term $\partial_t \bbz^n$ is approximated with BDF formulae of order $k\geq 1$, while the RHS of the ODE is solved (in this case) implicitly at time $t_n$. For a given $n \geq 1$ the BDF formulae read 
\begin{equation}\label{eq:bdf_k}
    \partial_t \bbz^n \approx \Xi^{\text{BDF}k}_n\left(\bbz\right) := \sum\limits_{p=0}^{k} \xi_p^{n,k} \bbz^{n-p},
\end{equation}
where the coefficients $\xi_p^{n,k}, p =0,\dots,k$ can be obtained by a Taylor expansion in $t^n$ and are reported in Table~\ref{tab:bdf_k} for $k=2,3$ for variable time-steps. 
\begin{table}
\caption{BDF$2$ and BDF$3$ expansion coefficients}   \label{tab:bdf_k}
\begin{center}
\bgroup
\begin{tabular}{|c|c|c|}\hline
        ~$p$~ & ~BDF$2$ coefficient~  & BDF$3$ coefficient \\    \hline 
        $0$ & $\frac{2\Delta_t^{n}+\Delta_t^{n-1}}{\Delta_t^{n}\left(\Delta_t^{n} + \Delta_t^{n} \right)}$  & $\frac{3 \left(\Delta_t^{n}\right)^2 + 4\Delta_t^{n}\Delta_t^{n-1}+2\Delta_t^{n}\Delta_t^{n-2} + \left( \Delta_t^{n-1}\right)^2 + \Delta_t^{n-1}\Delta_t^{n-2}}{\Delta_t^{n}\left( \Delta_t^{n}+\Delta_t^{n-1}\right)\left( \Delta_t^{n} + \Delta_t^{n-1}+\Delta_t^{n-2}\right)}$  \\ \hline
        $1$ & $- \frac{\Delta_t^{n} + \Delta_t^{n-1}}{\Delta_t^{n}\Delta_t^{n-1}}$ & $- \frac{\left( \Delta_t^{n}\right)^2+2\Delta_t^{n}\Delta_t^{n-1}+\Delta_t^{n}\Delta_t^{n-2}+\left( \Delta_t^{n-1}\right)^2+\Delta_t^{n-1}\Delta_t^{n-2}}{\Delta_t^{n}\Delta_t^{n-1}\left(\Delta_t^{n-1} + \Delta_t^{n-2}\right)}$\\ \hline
        $2$ & $\frac{\Delta_t^{n}}{\Delta_t^{n-1}\left( \Delta_t^{n} + \Delta_t^{n-1}\right)}$ & $\frac{\Delta_t^{n}\Delta_t^{n-1} +\Delta_t^{n}\Delta_t^{n-2} +\left(\Delta_t^{n}\right)^2 }{\Delta_t^{n-1}\Delta_t^{n-2}\left( \Delta_t^{n} + \Delta_t^{n-1}\right)}$\\ \hline
        $3$ & & $- \frac{\left(\Delta_t^{n}\right)^2+\Delta_t^{n}\Delta_t^{n-1}}{\Delta_t^{n-2}\left(\left(\Delta_t^{n-1}\right)^2+2\Delta_t^{n-1} \Delta_t^{n-2} +\Delta_t^{n}\Delta_t^{n-1} + \left(\Delta_t^{n-2}\right)^2 + \Delta_t^{n}\Delta_t^{n-2}\right) }$\\ \hline
\end{tabular}
\egroup
\end{center}
\end{table}

We now have all the necessary ingredients to introduce the time--adaptive algorithm based on the BDF$2$ and BDF$3$ schemes for the problem~\eqref{eq:general_nonlinear_problem}.

\subsection{Time-adaptive algorithm and temporal error estimators}

First of all, we assume that the continuous problem~\eqref{eq:general_nonlinear_problem} is well approximated by a mixed Finite Element (FE) method with appropriately chosen well-defined triangulations of the domains of interest and FE spaces $W_h \subset W$ and $W_{0, h} \subset W_0$ which are inf-sup stable to ensure the well-posedness of the FE formulation; we refer to~\cite{richter2017} for more details. In the following exposition, we will use $\bbz_h(t)$ to denote the FE approximation of the solution $\bbz(t)$ at time $t\in [0,T]$.

As mentioned in the introduction, time--adaptive algorithms are of extreme importance for complex physical models to produce the numerical simulations in feasible computational time. In this work, we resort to a methodology proposed already in 1961 in~\cite{ceschino1961modification}, reported in \cite{hairer1987solving}, and successfully studied in the context of FSI problems in~\cite{MAYR201855}. The core idea of the method lies in the next time-step size prediction: 
\begin{equation}\label{eq:step_size_prediciton}
    \Delta_t^\ast = \min \left\{ \Delta_t^{\max}, \max\left\{ \min\left\{ \kappa_{max}, \max\left\{ \kappa_{\min} , \kappa_s \kappa^\ast \right\} \right\} \Delta_t^n \right\}, \Delta_t^{\min} \right\}.
\end{equation}

The main component in~\eqref{eq:step_size_prediciton} is the scaling factor $\kappa^\ast$ which is defined as follows:
\begin{equation}\label{eq:kappa_star}
    \kappa^\ast = \left( \frac{\varepsilon}{\textit{est}_{n+1}}\right)^{\frac{1}{q+1}},
\end{equation}
where $q$ is the desired order of accuracy, $\varepsilon$ is the user-defined tolerance and $\textit{est}_{n+1}$ is the measure of the local temporal error.  Forgetting all other terms, the new time-step would be defined by $\Delta _t^* = \kappa^* \Delta _t^n$.
Then, one needs to add several safety factors. $\kappa_{min}$ and $\kappa_{\max}$ are user--specified parameters which indicate the minimal and the maximal ratio ${\Delta_t^*}/{\Delta_t^n}$ by which the time-step is allowed to increase or decrease. 
$\kappa_s<1$ keeps the local error incorporated in $\kappa^\ast$ away from the set tolerance. 
Furthermore, the bounds $\Delta_t^{\min}$ and $\Delta_t^{\max}$ are defined in order not to have an over/undershoot in the time--step prediction and they are usually chosen based on the physical knowledge of the problem at hand. 

We refer to~\cite{MAYR201855} for more insights on the general pipeline of the presented time--adaptive step selection, the choice of user--defined parameters and tolerances. 
Instead, we will focus on presenting a local temporal error estimator based on the interplay of implicit BDF$2$ and BDF$3$ schemes. 
Following the notation introduced in the equation~\eqref{eq:general_nonlinear_problem}, we denote by $\bbz_h^{n, \text{BDF$k$}}$, $k=2,3$ the solution to the following implicit in time problem:
\begin{equation}\label{eq:general_bdfk_implicit_scheme}
    R\left(\Xi_n^{\text{BDF}k}\left(\bbz_h^{\text{BDF$k$}} \right), \bbz_h^{n, \text{BDF$k$}}\right)=0 \quad \text{in} \ W_{0,h}', \quad n=3,4,\dots ,
\end{equation}
where the values $\bbz_h^{1, \text{BDF$k$}}$ and $\bbz_h^{2, \text{BDF$k$}}$ for $k=2,3$ are obtained by the implicit Euler and BDF$2$ schemes with constant time-step $\Delta_t^{\min}$, respectively, and the value $U_h^{0, \text{BDF$k$}}, k=2,3$ comes from the prescribed ICs.  As we aim to construct a second--order time advancing scheme, the value of parameter $q$ in~\eqref{eq:kappa_star} is chosen to be $q=2$.

The choice of using BDF implicit schemes is beneficial for many aspects for highly--nonlinear coupled problems, as the FSI problem~\eqref{eq:weak_time_continuous_fsi}. With respect to fully implicit multi--stage schemes, BDF implicit schemes lead to smaller nonlinear systems and they do not require the modification of the residual functional~\eqref{eq:general_nonlinear_problem} for different orders \cite{hairer1987solving}.
Moreover, both problems~\eqref{eq:weak_time_continuous_ns} and~\eqref{eq:weak_time_continuous_fsi} are of the form of differential--algebraic equations, i.e. some of the components (namely, the fluid pressure) or subequations (e.g., incompressibility constraints) do not contain the time derivative. 
The explicit schemes might not be able to capture fast transient dynamics strongly dominated by these components.   

For the above reasons, we propose the following form of the local temporal error estimator in~\eqref{eq:kappa_star}:
\begin{equation}\label{eq:estimator_implicit}
    \textit{est}_{n+1} = \max\limits_{u \in \bbz} ||u_h^{n+1,\text{BDF}2} - u_h^{n+1,\text{BDF}3}||_{L^2(\Omega)},
\end{equation}
that is the maximum absolute $L^2$--error over all the components of the problem~\eqref{eq:general_nonlinear_problem}. We will be referring to this estimator as \textit{implicit time estimator}.

The obvious issue of the estimator~\eqref{eq:estimator_implicit} is that it requires a lot of computational resources as it needs two non--linear highly--dimensional solvers. In order to mitigate this, we propose instead of solving for implicit BDF$3$ solution an approximation as a one--step Newton correction of the BDF$2$ solution as follows: find $\Tilde{\bbz}_h^{n+1,\text{BDF}3} = \bbz_h^{n+1,\text{BDF}2} + \delta \bbz_h^{n+1,\text{BDF}3}\in W$, where $\delta \bbz_h^{n+1,\text{BDF}3}\in W_0$ is the solution of the following linearised equation 
\begin{eqnarray}
 \nonumber J\left[ R\left(\Xi_{n+1}^{\text{BDF}3}\left(\bbz_h^{ \text{BDF$2$}} \right), \bbz_h^{n+1, \text{BDF$2$}}\right) \right]\left( \delta \bbz_h^{n+1,\text{BDF}3} \right) 
    \\   \label{eq:newton_bdf3_correction} = - R\left(\Xi_{n+1}^{\text{BDF}3}\left(\bbz_h^{ \text{BDF$2$}} \right), \bbz_h^{n+1, \text{BDF$2$}}\right) \quad \text{in } W_0'.
\end{eqnarray}
In the equation above $J$ represents the Jacobian operator which is computed by means of automatic differentiation and already used for the nonlinear solver of the BDF2. With this the time estimator~\eqref{eq:estimator_implicit} is approximated by the following:
\begin{equation}\label{eq:estimator_explicit}
    \textit{est}_{n+1} = \max\limits_{u \in \bbz} ||u_h^{n+1,\text{BDF}2} - \Tilde{ u}_h^{n+1,\text{BDF}3}||_{L^2(\Omega)},
\end{equation}
which requires only one nonlinear solver to get the BDF$2$ solution and one linear solver for~\eqref{eq:newton_bdf3_correction}. We will be referring to this estimator as \textit{linear-implicit (LI) time estimator}.

It is possible to prove, relying on approximation properties of BDF$k$ schemes, $k=2,3$ and \textit{a priori} estimates as the ones, for instance, in~\cite{PRUSAK2024253}, that the LI approximation defined by~\eqref{eq:newton_bdf3_correction} leads to the same order of accuracy as the implicit BDF$3$ solver~\eqref{eq:general_bdfk_implicit_scheme}. 
In the next section, we will show that estimators~\eqref{eq:estimator_implicit} and~\eqref{eq:estimator_explicit} have the same behaviour.

Finally, the next time-step prediction is defined as:
\begin{equation}\label{eq:next_time_step}
    \Delta_t^{n+1} := \alpha_0 \Delta_t^n + \alpha_1 \Delta_t^\ast,
\end{equation}
where the weights $\alpha_0,\alpha_1\geq 0$ satisfy $\alpha_0+\alpha_1=1$ and are chosen in order to prevent an overshooting effect which may lead to the need for significant reduction of the following time-steps. In the numerical results, we set $\alpha_0=0.3$ and $\alpha_1=0.7$.
As prescribed in \cite{MAYR201855}, we reperform the timestep evaluation until the condition $\textit{est}_{n+1} < \varepsilon$ is satisfied and, in any case, for a maximum of 5 iterations.
\section{Numerical simulations}\label{sec:simulations}

In this section, we will provide a few numerical test cases to show the efficiency of the presented time--adaptive algorithm. We will provide the numerical results on two benchmark problems: the backward--facing step flow and a two--dimensional haemodynamic FSI test case. 
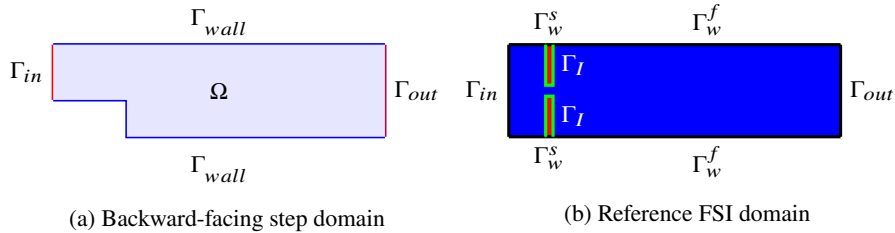
\begin{figure}
    \centering
\begin{subfigure}[h] {0.49\textwidth}
\begin{tikzpicture}[scale=0.49]
    \draw[draw=red,very thick] (0,1) -- (0, 1.75) node [anchor=east]{$\Gamma_{in}$}  -- (0,2.5);
    \draw[draw=blue,very thick] (0,2.5) --  (4.5, 2.5) node [anchor=south] {$\Gamma_{wall}^{\vphantom{f}}$}-- (9,2.5); 
    \draw[draw=purple,very thick](9,2.5) -- (9, 1.25) node [anchor=west] {$\Gamma_{out}$}-- (9,0);
    \draw[draw=blue,very thick] (9,0) -- (4.5, 0) node [anchor=north]{$\Gamma_{wall}^{\vphantom{f}}$} -- (2,0)--(2,1)--(0,1);
    \fill[fill=blue!10] (0,1)  -- (0,2.5) -- (9,2.5) -- (9, 1.25) -- (9,0)  -- (2,0)--(2,1)--(0,1);    
    \draw (4.5, 1.25) node {$\Omega$};
\end{tikzpicture}
\caption{Backward-facing step domain\label{fig:domain_bfs}
\captionsetup{justification=centering}
}
\end{subfigure}    
\begin{subfigure}[h] {0.5\textwidth}
\begin{tikzpicture}[scale=0.49]
    \fill[fill=blue] (0,0)  -- (1,0) -- (1,1.1) -- (1.2,1.1) -- (1.2,0)--(9,0)--(9,2.5) --(1.2, 2.5)-- (1.2, 1.4) -- (1,1.4) -- (1, 2.5) -- (0, 2.5) -- (0,0); 
    \fill[fill=red] (1,0) -- (1,1.1)-- (1.2,1.1) -- (1.2, 0) -- (1,0) ;
    \fill[fill=red] (1,2.5) -- (1,1.4)-- (1.2,1.4) -- (1.2, 2.5) -- (1,2.5) ;
    \draw[draw,very thick] (0,0) -- (0, 1.25) node [anchor=east]{$\Gamma_{in}$}  -- (0,2.5);
    \draw[draw,very thick] (1.2,2.5) --  (5.4, 2.5) node [anchor=south] {$\Gamma_{w}^f$}-- (9,2.5);
    \draw[draw,very thick] (1.2,0) --  (5.4, 0) node [anchor=north] {$\Gamma_{w}^f$}-- (9,0);
    \draw[draw,very thick](9,2.5) -- (9, 1.25) node [anchor=west] {$\Gamma_{out}$}-- (9,0);
    \draw[draw,very thick] (1,2.5) --  (1.1, 2.5) node [anchor=south] {$\Gamma_{w}^s$}-- (1.2,2.5);
    \draw[draw,very thick] (1,0) --  (1.1, 0) node [anchor=north] {$\Gamma_{w}^s$}-- (1.2,0); 
       \draw[draw, very thick](0,0)--(1, 0);
    \draw[draw, very thick](0,2.5)--(1, 2.5);
     \draw[draw=green,very thick] (1,0) --  (1, 1.1)  --(1.2,1.1) -- (1.2,0) ; 
    \draw[draw=green,very thick] (1,2.5) --  (1, 1.4)  --(1.2,1.4) -- (1.2,2.5) ; 
\node[anchor=west, color=white] at  (1.2, 0.55) {$\Gamma_{I}$};
\node[anchor=west, color=white] at  (1.2, 1.95) {$\Gamma_{I}$};
\end{tikzpicture}
\caption{Reference FSI domain } \label{fig:domain_fsi}%In blue: the reference fluid domain $\Omega^f$. In red: the reference solid domain $\Omega^s$. In green: the fluid--structure interface $\Gamma_I$.}
\captionsetup{justification=centering}
\label{fig:domain_valves_fsi}
\end{subfigure}   
\caption{Domains of interest for the CFD (a) and FSI (b) test cases\label{fig:domains}
}
\end{figure}

\textbf{Backward-facing step flow CFD test case.} Figure~\ref{fig:domain_bfs} represents the physical domain of interest. The upper part of the channel has a length of $18 cm$, the lower part $14 cm$; the height of the left chamber is $3 cm$, and the height of the right one is $5 cm$. The spatial discretisation is carried out by FE method using Taylor-Hood $\mathbb{P}_2-\mathbb{P}_1$ FE pair with $27, 890$ degrees of freedom (DoFs). We impose homogeneous Dirichlet boundary conditions (BCs) on the top and the bottom walls of the boundary $\Gamma_{wall}$ for the fluid velocity, the homogeneous Neumann (free--outflow) conditions on the outlet $\Gamma_{out}$  on this portion of the boundary, and the parabolic profile $u_{in}(x,y,t) = \left( \varphi(t) \frac{20}{9} (y-2)(5-y),0  \right)^T$ on $\Gamma_{in}=\lbrace (x,y): x= 0, y \in [2,5]\rbrace, \ t=[0,2]$, where $\varphi\in \mathcal C^1(\mathbb R^+)$ is given by
\begin{equation*}
    \varphi(t) = \begin{cases}
        \frac{1}{2}\left( 1- \cos(\pi t)\right) & \text{for }t\leq 1, \\
        1 & \text{for }t >1.
    \end{cases}
\end{equation*}
Below we will provide the numerical simulation of the proposed time--adaptive algorithm for the value of the viscosity parameter $\nu=0.05$ resulting in Reynolds numbers equal to $300$. We will refer to this test case as \textbf{CFD-300} in the following.

\textbf{Two--dimensional haemodynamics FSI test case (\textbf{FSI-H2})}. \begin{figure}
    \centering
    \begin{subfigure}{0.49\textwidth}
        \includegraphics[width=\textwidth]{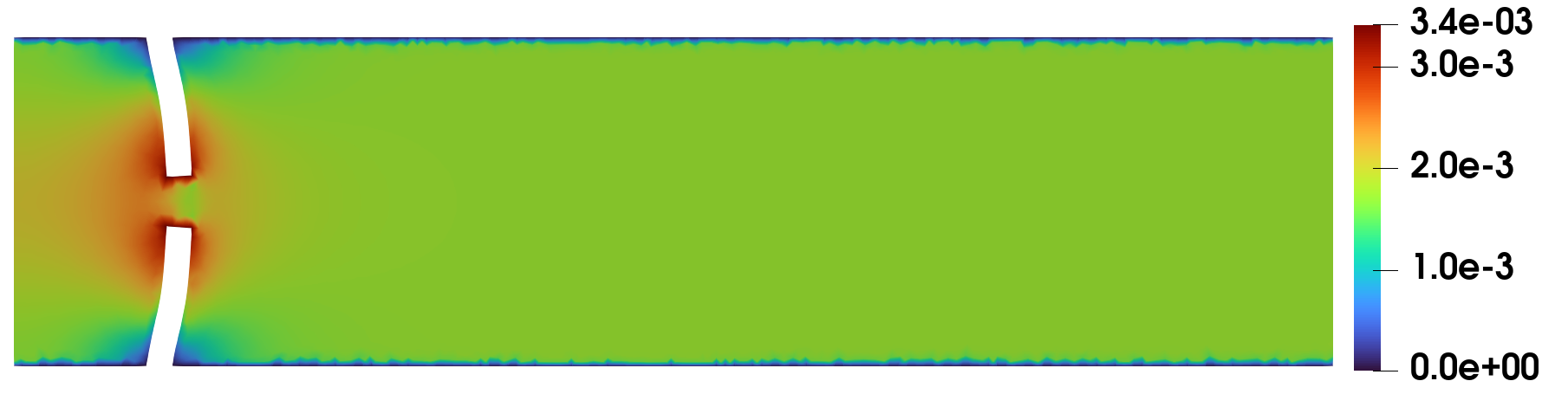}
        \caption{Fluid velocity}
    \end{subfigure}
    \begin{subfigure}{0.49\textwidth}
        \includegraphics[width=\textwidth]{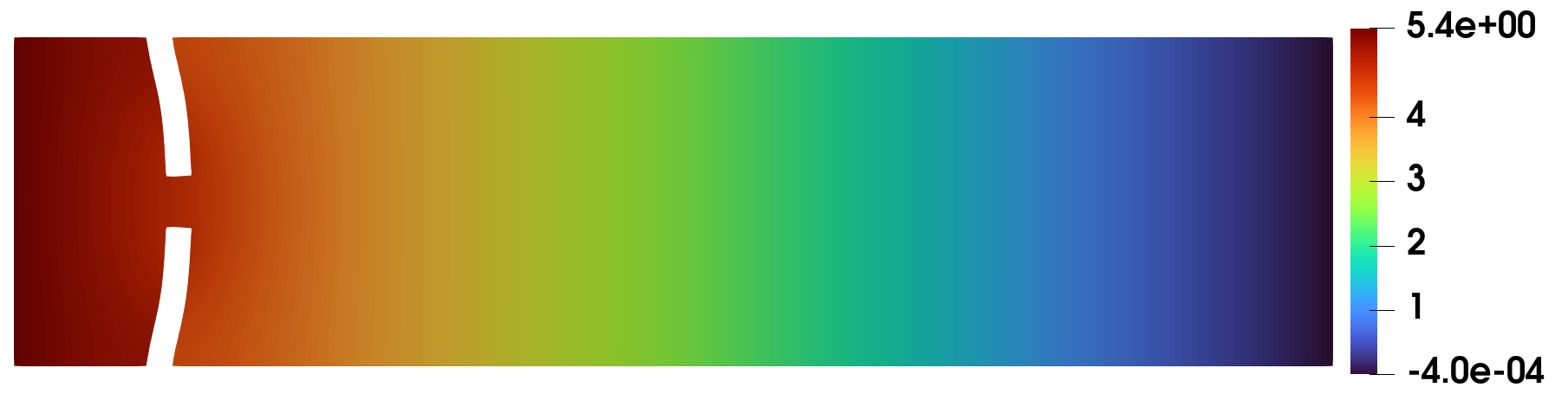}
        \caption{Fluid pressure}
    \end{subfigure}
    % \begin{subfigure}{0.19\textwidth}
    %     \includegraphics[scale=0.01]{Figures/d_s_rom_100.png}
    % \end{subfigure}
    \caption{\textbf{FSI-H2}  solution at time instance $t=0.01$ in the current fluid configuration}
    \label{fig:fsi_solution_001}
\end{figure}
Figure~\ref{fig:domain_fsi} represents the physical domain in the reference configuration: the reference fluid domain $\Omega^f$ in blue, the reference solid domain $\Omega^s$ in red and the fluid--structure interface $\Gamma_I$ in green. the fluid domain is $2.5 cm$ in height $10 cm$ long; the leaflets are situated $1 cm$ downstream the inlet boundary $\Gamma_{in}$, they are $0.2 cm$ thick and $1.1 cm$ in height. The values of physical parameters are the following: $\nu_f = 0.035, \rho_f = 1, \rho_s = 1.1, \lambda_s = 8\cdot 10^{5}$ and $\mu_s = 4 \cdot 10^5$. The spatial discretisation is carried out by FE method using a generalised Taylor-Hood~\cite{richter2017} FE triple $\mathbb{P}_2-\mathbb{P}_1-\mathbb{P}_2$ with $67, 390$ DoFs.
We consider zero ICs for fluid velocity and the structure displacement, homogeneous Dirichlet BCs on $\Gamma^f_{w}$ for the fluid velocity and on $\Gamma_w^s$ for the structure displacement, and homogeneous Neumann (free--outflow) conditions on $\Gamma_{out}$. A pressure impulse $\bbu_N^f(x,t) = -p_{in}(t)\textbf{n}_f(x) = \left(p_{in}(t), 0\right)^T, \forall x\in \Gamma_{in}, \forall t \in [0,2]$, is applied as a Neumann condition at $\Gamma_{in}$ with: 
\begin{equation*}
p_{in}(t) = 
\begin{cases}
        5\left( 1- \cos(\frac{\pi t}{0.2})\right) & \text{for }t\leq 0.1, \\
        5 & \text{for }t >0.1.
\end{cases}
\end{equation*}

In Fig.~{\ref{fig:fsi_solution_001}}, we provide a plot of the fluid subcomponents of the \textbf{FSI-H2} test case at the time $t=0.01$ in the current fluid domain configuration. The magnitude of the displacement at this time instance is of order $10^{-6}$, and it can be seen that the pressure field dominates the dynamics of the problem. 
\begin{figure}[h]
    \centering
    \begin{subfigure}{0.49\textwidth}
     \includegraphics[width=\textwidth]{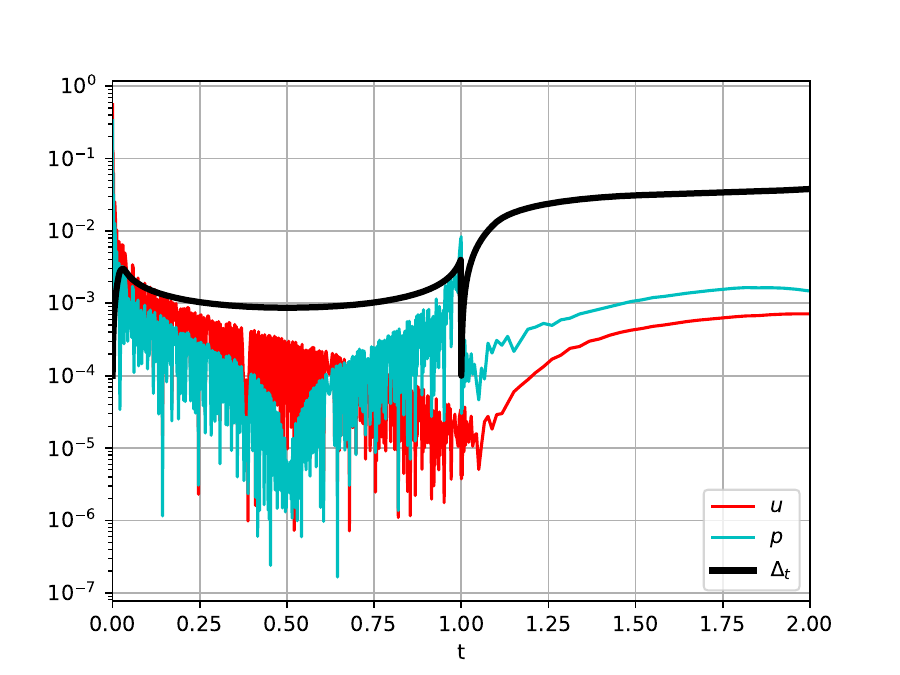}
    \caption{\textbf{CFD-300}: errors and timestep}\label{fig:adaptive-NS-err}
    \end{subfigure}
    % \label{fig:enter-label}
    \begin{subfigure}{0.49\textwidth}
    \includegraphics[width=\textwidth]{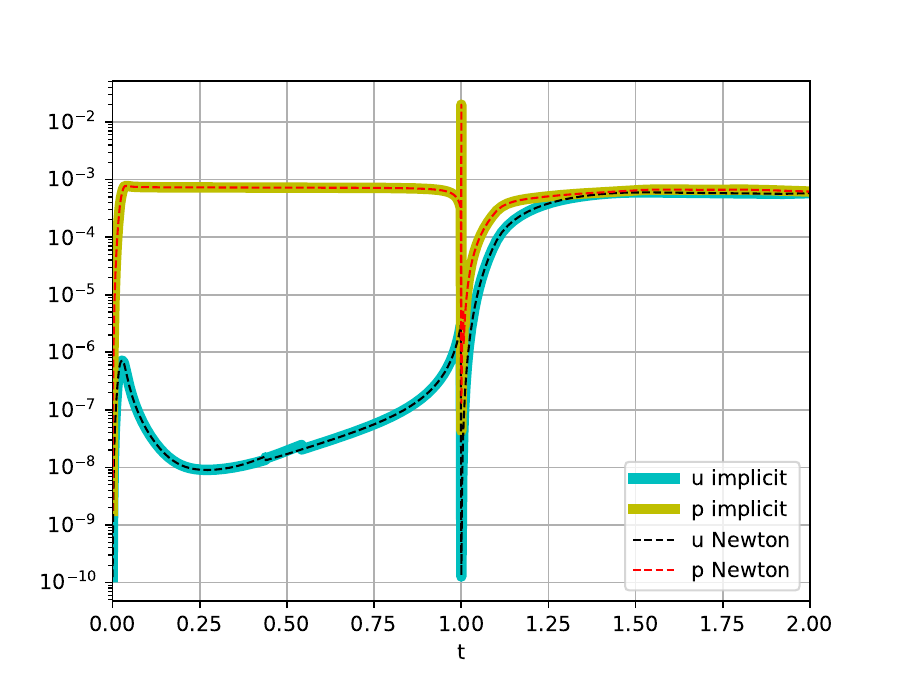}
    \caption{\textbf{CFD-300}: error estimators}\label{fig:adaptive-NS-errest}
    \end{subfigure}
        \begin{subfigure}{0.49\textwidth}
     \includegraphics[width=\textwidth]{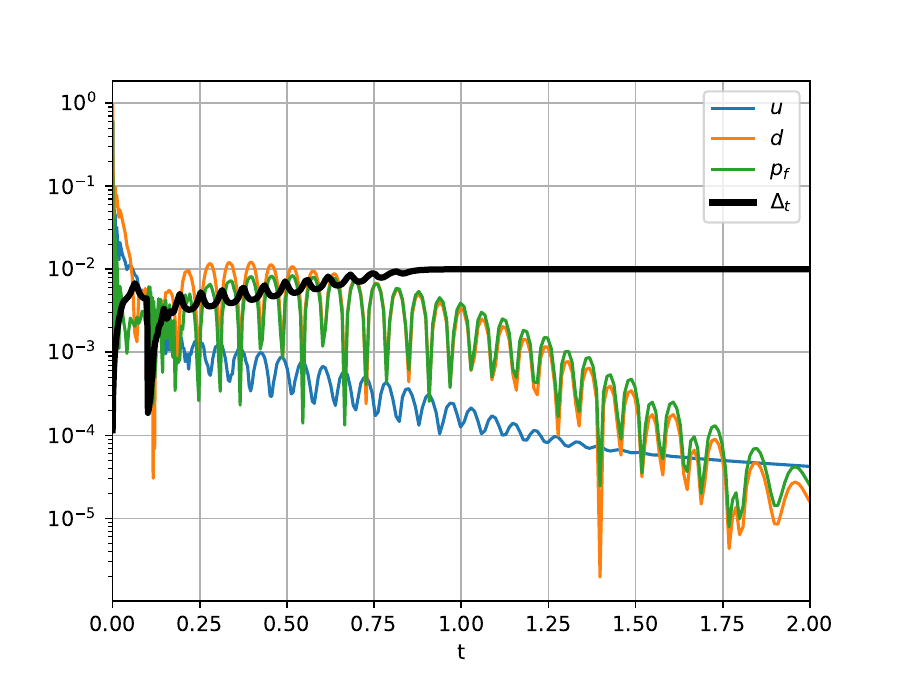}
    \caption{\textbf{FSI-H2}: errors and timestep}\label{fig:adaptive-FSI-err}
    \end{subfigure}
    % \label{fig:enter-label}
    \begin{subfigure}{0.49\textwidth}
    \includegraphics[width=\textwidth]{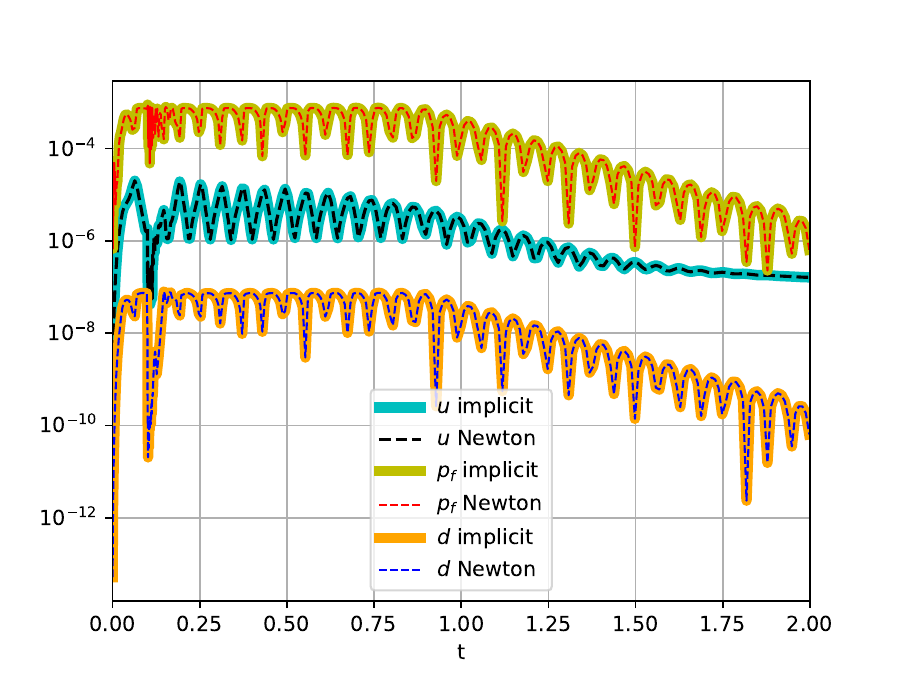}
    \caption{\textbf{FSI-H2}: error estimators}\label{fig:adaptive-FSI-errest}
    \end{subfigure}
    \caption{Adaptive time-steps distribution and relative errors w.r.t. constant-timestep solution (left) and the comparison of implicit and LI time estimators (right)}\label{fig:adaptive}
\end{figure}
\begin{table}[h!]
\caption{Computational cost comparison between constant BDF2 and time adaptive algorithms: total simulation time (left) and time of one evaluation of the error estimator (mean $\pm$ std) }   \label{tab:comp_times} 
\begin{center}
\bgroup
\begin{tabular}{|c||c|c||c|c||c|c|c|}\hline
         & \multicolumn{2}{c||}{Computational time} & \multicolumn{2}{c||}{Estimator cost} & \multicolumn{3}{c|}{Steps number}  \\    \hline 
        ~Test case~ & ~Constant~& ~LI adaptive~ & ~Implicit~&~LI~ & ~Const.~ & ~Impl.~ & ~LI~ \\    \hline 
        \textbf{CFD-300} & 12 hours & 2 hours & 2.1  $\pm$ 0.2 sec  & 1.68 $\pm$ 0.1 sec  &~20,000~ &~976~&~976~  \\ \hline 
        \textbf{FSI-H2} & 132 hours & 5 hours & 37 $\pm$ 4 sec &  25  $\pm$ 2 sec  &~20,000~ &~347~&~360~ \\ \hline
\end{tabular}
\egroup
\end{center}
\end{table}

\textbf{Validation of the time--adaptive algorithm}. Having established the test cases, we will now provide a numerical analysis of the time--adaptive algorithm described in Section~\ref{sec:time_estimator}. In both test cases we fix the following parameters entering the timestep prediction~\eqref{eq:step_size_prediciton}: $\kappa_{\max} = 1.5, \ \kappa_{\min} = 0.1, \ \kappa_s = 0.9$ and the tolerance $\varepsilon=10^{-3}$ in~\eqref{eq:kappa_star}. For \textbf{CFD-300} we choose $\Delta_t^{\min} = 10^{-4}$ and $\Delta_t^{\max} = 10^{-1}$ whereas for  \textbf{FSI-H2} we set $\Delta_t^{\min} = 10^{-4}$ and $\Delta_t^{\max} = 10^{-2}$. 
In Figure~\ref{fig:adaptive-NS-err} for \textbf{CFD-300}, we observe the timestep chosen by the algorithm and the error w.r.t. a reference solution (computed via $\Delta_t^{\min}$ constant--timestep BDF2 scheme). We observe that the global error is somewhat under control during the whole simulation. The errors and the time steps are particularly sensitive to the discontinuity of the second derivative in time of the BCs at $t=1$. The error estimator in Figure~\ref{fig:adaptive-NS-errest} is indeed based on Taylor expansions of the second order and it correctly detects a singularity and decreases the timestep in that point, maintaining the error of the pressure (that varies widely at that time) under control. 
The differences between the fully implicit and the LI error estimators are negligible in terms of accuracy, but they sum up to a strong reduction in computational costs as evidenced in Table~\ref{tab:comp_times}. 

% \begin{figure} 
%         \begin{subfigure}{0.49\textwidth}
%      \includegraphics[width=\textwidth]{Figures/Timesteps_relerr_fsi.pdf}
%     \end{subfigure}
%     % \label{fig:enter-label}
%     \begin{subfigure}{0.49\textwidth}
%     \includegraphics[width=\textwidth]{Figures/Estimators_fsi.pdf}
%     \end{subfigure}
%     % \caption{Caption}
%     \caption{\textbf{FSI-H2} test case: adaptive time-steps distribution and relative errors w.r.t. constant-timestep solution (on the left) and the comparison of implicit and Newton-explicit time estimators (on the right)}\label{fig:adaptive-FSI}
% \end{figure}

Similarly, for \textbf{FSI-H2} in Figure~\ref{fig:adaptive-FSI-err} we notice that the error is always under control. In this simulation, the relative errors are comparable among components in particular for pressure and displacement. On the other hand, the error estimator in Figure~\ref{fig:adaptive-FSI-errest} is led by the pressure and it is of paramount importance to include such a component in the error estimator. 
Again, there is essentially no difference between the fully implicit estimator and the LI one, but great computational saving is achieved as can be seen in Table~\ref{tab:comp_times}.
%

%
%Vale mettere footnotesize? O è cheating? Magari non cosi´ piccolo pero un smallsize sicuramente  ++++. come cazz si faceva haha Io ho messo footnotesize non è così piccolo danke =) una referenza sborda PD uno scriptsize? Cosi´ sta tutto. Fatemi sapere se e´ troppo piccolo.. D: Un po' estremo scriptsize 
%M: eh lo so ma come facciamo altrimenti? Non mi vengono in mente alternative a meno di non abbreviare qualcosa e mangiarci 3 righe..
% D: vado a tagliare cose M: vai mago del taglia&cuci <3
\footnotesize
% \scriptsize
\subsection*{Acknowledgements} 
This work was supported by the European Union’s Horizon 2020 research and innovation programme under the Marie Sklodowska–Curie Actions [872442], PNRR NGE iNEST project and the European High-Performance Computing Joint Undertaking (JU) [955558]. IP has been funded by a SISSA fellowship within the projects TRIM [G95F21001070005]. DT was funded by a SISSA Mathematical fellowship and by the Ateneo Sapienza projects 2022 “Approssimazione numerica di modelli differenziali e applicazioni” and 2023 “Modeling, numerical treatment of hyperbolic equations and optimal control problems”.
DT and GR are members of the INdAM Research National Group of Scientific Computing (INdAM-GNCS). This research was funded in part by the Austrian Science Fund (FWF) project 10.55776/ESP519 (MN). For open access purposes, the author has applied a CC BY public copyright license to any author-accepted manuscript version arising from this submission.

\normalsize
\bibliographystyle{siam}
\bibliography{main}

\end{document}